\documentclass[11pt]{amsart}

\usepackage[T1]{fontenc}
\usepackage[utf8]{inputenc}
\usepackage{lmodern}
\usepackage{microtype}
\usepackage{amsmath,amssymb,amsfonts,mathtools}
\usepackage{enumitem}
\usepackage[colorlinks=true,linkcolor=blue,citecolor=blue,urlcolor=blue]{hyperref}
\usepackage{geometry}
\usepackage{mathrsfs}

\newtheorem{theorem}{Theorem}[section]
\newtheorem{proposition}[theorem]{Proposition}
\newtheorem{lemma}[theorem]{Lemma}
\newtheorem{corollary}[theorem]{Corollary}
\newtheorem{definition}[theorem]{Definition}
\newtheorem{remark}[theorem]{Remark}

\newcommand{\R}{\mathbb R}
\newcommand{\tr}{\operatorname{tr}}
\newcommand{\Lip}{\operatorname{Lip}}
\newcommand{\diag}{\operatorname{diag}}
\newcommand{\Sym}{\operatorname{Sym}}
\newcommand{\Skew}{\operatorname{Skew}}

\newcommand{\LL}{\mathrm{L}}
\title{Rank-one convex Hessians and the extension problem}
\author{Gabriele Cassese}
\date{}

\begin{document}

\begin{abstract}
We prove that the Hessian of a rank-one convex integrand $f:\R^{d\times m}\to \R$ is measure-valued. We also 
show that rank-one convex and quasiconvex integrands defined on subspaces such as $\Sym(2)$ do not 
necessarily admit extensions to the full ambient space.
\end{abstract}
\maketitle

\section{Introduction}
Rank-one convex integrands have long played an important role in vector-valued calculus of variations.
Despite this, not much is known about their structure. In particular, it is easy to prove that
they are locally Lipschitz (provided the domain is spanned by its rank-one cone), but the parallel with convexity suggests they should have
measure-valued Hessians and a.e. defined pointwise Hessians. More precisely, given a
subspace $L$ of $\R^{d\times m}$, one can ask
\begin{center}
      $(Q_1)$\quad \textit{Does a rank-one convex integrand $f:L\to \R$ have a measure-valued Hessian?}
\end{center}
\begin{center}
      $(Q_2)$\quad \textit{Does a rank-one convex integrand $f:L\to \R$ have an a.e. pointwise Hessian?}
\end{center}
\begin{center}
$(Q_3)$\textit{Does a (strictly) rank-one convex integrand 
$f:L\to \R$ admit a rank-one convex extension $\hat f:\R^{d\times m}\to \R$?}
\end{center}
$Q_1$ was originally posed, and answered negatively for $L=\Sym(2)$, by Conti, Faraco, Maggi and M\"uller in \cite{Faraco}. This result was later generalised to $\Sym(n\ge 2)$
by Kirchheim and Kristensen in \cite{KK}; more generally their Theorem 1.4 gives necessary and sufficient conditions
on $L$ for $Q_1$ to have a positive answer, provided $L$ satisfies a certain assumption. Unfortunately,
$\R^{m\times d}$ does not satisfy this assumption.

$Q_2$ was recently answered positively by Hirsch
in \cite{Hirsch}, proving that rank-one convex integrands $f:L\to \R$ are, in a suitable sense,
second-differentiable almost everywhere provided the rank-one cone of $L$ spans $L$. In particular,
$\Sym(2)$ gives a positive answer to $Q_2$ but not to $Q_1$. 

$Q_3$ also originates from \cite{Faraco}, and has attracted
additional attention in its restricted form applied to 
quasiconvex integrands (\cite{DalMaso}, \cite{Cagnetti}).
Much is known under additional regularity hypotheses on
$f$,
but the question has remained unanswered in its full generality, particularly when $D^2f$ is locally unbounded.
\par
Nonetheless, the question of whether $Q_1$ holds if $L$ is the whole matrix space, say $\R^{2\times2}$, has remained
open; there is a natural relation with $Q_3$:
 if an extension procedure does exist, say from $\Sym(2)$
 to $\R^{2\times 2}$, then one could hope to explore 
 $Q_1$ by analysing the behaviour of the extension of integrands on $\Sym(2)$ with locally unbounded Hessians.\par
 In this paper we answer the question by proving that rank-one convex integrands
on the whole matrix space have measure-valued Hessians.

\begin{theorem}\label{thm:hess-meas}
Let $f:\R^{d\times m}\to\R$ be a rank-one convex integrand. Its
distributional Hessian is a matrix-valued Radon measure. More precisely,
\[
 |D^2f|(B_r(A))\lesssim_{d,m}
 r^{dm-1}\Lip(f;B_{8r}(A)).
\]
\end{theorem}

Combined with differentiation of measures and the local compactness of
separately convex functions, this gives an alternative proof of the result of Hirsch \cite{Hirsch}, in the case
of integrands defined on the full ambient space.

\begin{corollary}\label{cor:hess-point}
Let $f:\R^{d\times m}\to\R$ be a rank-one convex integrand, and write
$D^2f=H_f\mathcal L^{dm}+(D^2f)^s$. Then, for almost every $A$,
\[
 f(A+Z)=f(A)+Df(A)(Z)+\frac12 H_f(A)[Z,Z]+o(|Z|^2).
\]
\end{corollary}

On the other hand, we prove that $Q_3$ has a negative answer in all but the trivial cases. More precisely, say $L$ has the rank-one extension property 
if and only if each rank-one convex integrand $f:L\to \R$ admits an extension. Then
\begin{theorem}\label{thm:ext}
      $L$ has the rank-one extension property if and only if $L$ is a rectangular subspace, i.e. of the form 
      $U\otimes V$ where $U$ is a subspace of $\R^{m}$ and $V$ of $\R^d$. If this is not the case, 
      there exists a uniformly strictly rank-one convex integrand $f:L\to \R$ which does not admit
      a rank-one convex extension to the ambient space.
      If $L\subset \R^{2\times 2}$, one may even assume 
      $f$ to be strictly uniformly quasiconvex.
\end{theorem}

The Hessian result can be generalised as follows. Given a symmetric quadratic
form $q$ on $\R^N$, set
\[
 \mathcal D_q\mathrel{:=}\{v\in\R^N:q(v,v)=0\}.
\]
An integrand is $\mathcal D_q$-convex if its restriction to every line
in a direction belonging to $\mathcal D_q$ is convex.

\begin{proposition}\label{prop:classification}
Let $q$ have positive and negative indices $p$ and $n$, respectively,
and assume $n>0$. Then:
\begin{enumerate}
 \item If $\min(p,n)\le1$, there exists a $C^1$, $\mathcal D_q$-convex
 integrand whose distributional Hessian is not measure-valued.
 \item If $\min(p,n)\ge2$, every finite $\mathcal D_q$-convex integrand
 has a measure-valued Hessian.
\end{enumerate}
\end{proposition}

Thus, for quadratic null cones, we complete the classification initiated
by Kirchheim--Kristensen \cite{KK}.

In sections \ref{sec:dc} and \ref{sec:ext} we turn our attention to two additional 
natural questions: whether suitably regular rank-one convex integrands can be locally written as the difference of 
two convex functions, see theorem \ref{thm:noncvx}, and to the problem of obtaining a rank-one convex extension for integrands defined on $\Sym(2)$ and more generally subspaces $L\subset \R^{m\times d}$ to the whole space, see \ref{thm:ext}.

\subsection{Notation and conventions}
Throughout, matrix norms are assumed to be Frobenius norms, 
with the inner product denoted by $A:B$ and defined as $A:B=\mathrm{Tr}(AB^{\top})$, and $B_r$ is the ball of radius $r$ centred at
zero unless a centre is specified. Constants may change from line to
line. We identify $\R^{2\times2}$ orthogonally with $\R^4$ so that
$\det A=\frac12x^\top Qx$, where $Q= \diag(1,1,-1,-1)$. In these coordinates, rank-one directions
are precisely the $Q$-null directions. Throughout the paper we will refer to continuous maps $f:L\subset \R^{m\times d}\to \R$
as integrands.\par 
We will make use of the following shorthand notation 
for sets of square matrices.
$\mathrm{Sym}(n)$ denotes the space of symmetric $n\times n$
 matrices and $\Sym A=1/2(A+A^{\top})$, $\mathrm{CO}^+(n)$ the set of 
$n\times n$ matrices defined as $\{\rho U:\rho\ge 0,\ U\in \mathrm{SO}(n)\}$
,
with $\mathrm{CO}^-(n)\mathrel{:=}\{\rho U:\rho\ge 0,\ U\in \mathrm{O}(n)\setminus \mathrm{SO}(n)\}$.
For $n=2$, $\mathrm{CO}^-(2)=(\mathrm{CO}^+(2))^\perp$. $\mathrm{Tri}(n)$ and $\mathrm{diag}(n)$
will denote, respectively, the space of $n\times n$
upper-triangular matrices and that of $n\times n$ diagonal matrices, and $\mathfrak{sl}(n)$ 
the space of trace-free $n\times n$ matrices.\par
Given a measure $\mu$, its barycentre will be denoted
$[\mu]$, and the integration pairing 
$\langle f,\mu\rangle=\int f\mathrm{d}\mu$.
\section{An overview of the problem and its solution}
Let us start by recalling the definition of a rank-one convex integrand.
\begin{definition}
Given a subspace $L$ of $\R^{m\times d}$, an integrand $f:
L\to \R$ is said to be rank-one convex if,
for all matrices $A,B\in L$ such that $\mathrm{rank}(A-B)= 1$, the map $t\mapsto f(tA+(1-t)B)$ is convex on $[0,1]$. $f$ is said to be strictly 
rank-one convex if $t\mapsto f(tA+(1-t)B)$ is strictly convex on $[0,1]$ and uniformly strictly 
rank-one convex if there exists $\varepsilon>0$
such that $f-\varepsilon|\cdot|^2$ is rank-one convex.
\end{definition}
More generally, given a closed cone 
$\mathcal{D}$ in $\R^N$ we define $\mathcal D$-convexity as
\begin{definition}
    A function $f:\R^N\to \R$ is said to be $\mathcal D$-convex if for any $v\in D$, $a\in \R^N$
    the function $t\mapsto f(a+tv)$ is convex on 
    $[0,1]$. It is said to be $\mathcal D$-affine
    if both $f$ and $-f$ are $\mathcal D$-convex.
\end{definition}
The usefulness and the applications of this concept are manifold: we refer the reader to \cite{Rindler} for a thorough exposition of the subject.\par
Consider for simplicity's sake
a rank-one convex integrand $f$ on $\R^{2\times 2}$.
Rank-one convexity implies that $f$ satisfies the distributional Legendre-Hadamard condition:
for any rank-one matrix $H\in \R^{2\times 2}$
\begin{equation}
D^2 f(\cdot)[H,H]\ge 0\label{eq:LH}.
\end{equation}
Let us consider the standard basis of $\R^{2\times 2}$: $e_{11},e_{12},e_{21}, e_{22}$.
Applying \eqref{eq:LH} to the basis vectors
it follows that the diagonal components of the Hessian
$D^2f$, namely $f_{e_{11}e_{11}}, f_{e_{12}e_{12}}, \dots $
are nonnegative distributions, and hence Radon measures.
Taking other rank-one matrices in the Legendre-Hadamard condition it also follows that mixed derivatives sharing a column ($f_{e_{11}e_{21}}, f_{e_{12}e_{22}})$ can be written as differences of positive distributions and thus are Radon measures. For example
\[
4f_{e_{11}e_{12}}=D^2f[e_{11}+e_{12}, e_{11}+e_{12}]-D^2f[e_{11}-e_{12}, e_{11}-e_{12}].
\]

Thus the only obstacle to $D^2f\in \mathcal M_{loc}(\R^{2\times 2})$ are $f_{e_{11}e_{22}}$ and $f_{e_{12}e_{21}}$. Applying \eqref{eq:LH}
with $H=e_{11}+e_{12}+e_{21}+e_{22}$
yields
\[
f_{e_{12}e_{21}}+f_{e_{11}e_{22}}=\frac12D^2f[H,H]-\frac12\Delta f-
(f_{e_{11}e_{12}}+f_{e_{11}e_{21}}+f_{e_{12}e_{22}}+f_{e_{21}e_{22}})\]
and thus their sum is a measure.

The difficulty is therefore concentrated in a single component of the Hessian, which is invisible to
the Legendre-Hadamard inequalities and to the Laplacian, namely $f_{e_{11}e_{22}}-f_{e_{12}e_{21}}$. Note that
this reduction to a single coordinate should not be taken as a strong indicator of the truth value of $Q_1$, since
a similar result can be proved in the case of $\Sym(2)$
considered in \cite{Faraco}: namely, if we consider a rank-one convex integrand on $\Sym(2)$, of the form
\[
f\left(\begin{pmatrix}
      a&b\\ b&d
\end{pmatrix}\right)
\]
then by using suitable rank-one matrices in \eqref{eq:LH} it follows that $f_{aa}, f_{dd}, f_{ab}, f_{bd}$
and $f_{bb}+2f_{ad}$ are Radon measures, so the failure of $D^2 f\in\mathcal M_{loc}(\Sym(2))$ rests entirely
on the failure of $f_{bb}-2f_{ad}$ to be a measure.\par
The problem is clearly one of regularity theory for differential inclusions: given $f$ such that $D^2f\in \Gamma$
in the appropriate weak sense, does it imply $Df\in \mathrm{BV}_{\mathrm{loc}}$.
The ideas of this paper were deeply inspired by the notion of $k$-convexity and in particular by
the work \cite{TW}. In the study of $k$-convexity, one is interested in functions
$f:\R^{N}\to \R$ such that $\sigma_1(D^2f)\ge0,\dots, \sigma_k(D^2f)\ge 0$ in the appropriate weak sense, where
$\sigma_j$ denotes the $j$-th symmetric polynomial on the eigenvalues. Clearly, rank-one convex functions are
$1$-convex, i.e. subharmonic. Moreover, $2$-convexity is enough to ensure the Hessian is measure-valued (see \cite[Theorem 2.4]{ChaudhuriTrudinger} and \cite[Lemma 2.3]{TW}); informally, for a smooth function \footnote{The full argument is slightly more complicated as it requires dealing with mollifications, see the references above for a complete proof}
$\sigma_1(D^2f)\ge 0$ and $\sigma_2(D^2f)\ge 0$ imply
\[
|D^2f|\le \Delta f
\]
and since $\Delta f\ge 0$, it is a Radon measure.
\par
The use of $k$-convexity in the calculus of variations is not new to this paper, see \cite{FaracoZhong}.
On the other hand, our construction relies on a modified version of $k$-convexity, which we now briefly 
describe.\par
Our results rest on the observation that there exists (on $\R^{2\times 2}$) a positive\footnote{Meaning that it maps positive semidefinite matrices to nonnegative numbers} quadratic null-lagrangian polynomial $K$ such that for any rank-one convex integrand
$f$ we have $K(D^2f)\ge 0$. More precisely, one has (see Lemma \ref{lemma:cone}) the following decomposition of the Hessian:
\[
 D^2f=P+uQ,\qquad P\ge0,\qquad Q= \diag(1,1,-1,-1), \qquad u\in\R,
\]
and our polynomial $K$ will be positive and invariant under 
$u$-translations, i.e. $K(H+uQ)=K(H)$.
Since $|P|\le\tr P=\Delta f$, it remains to control $u$.
The null-lagrangian identity of $K$ gives an upper bound on 
$\|K(D^2 f)\|_{\LL^1}$, but this by itself is not sufficient to control
$u$ because of the translation invariance of $K$.
We apply the integral estimate for $K$ after quadratic
sup-convolution, motivated by Hirsch's application of the sliding paraboloid method to rank-one convex integrands \cite{Hirsch}. Where a paraboloid supports $f$ from above,
the construction of $K$ yields $K(D^2f_{\varepsilon})\gtrsim \varepsilon^2 u^4$ and the integral upper bound gives us the desired control.
A standard estimate for the excluded set will then yield
an integrable tail. \par
The result can be thought of as a regularity result for a differential inclusion:  $D^2f\in \Gamma$ (where $\Gamma$ denotes the quadratic forms that are positive-semidefinite on the rank-one set) implies $\nabla f\in \mathrm{BV}_{\mathrm{loc}}$. As a regularity result, it is rather peculiar
in its conclusion: most results concerning similar cone differential inclusions concern $L^p$ estimates
(see \cite{Boman1967,Boman1982,GuerraRaitaSchrecker}), and their techniques do not seem to easily extend to this case. Some results concerning $\mathrm{BV}$-regularity are known, see \cite{AlbertiAmbrosio}, particularly Theorem 5.3 and Remarks 5.4-5.5.

We briefly sketch some additional applications of the technique we implement, but we expect this approach (or rather, trick) to be useful to more
general regularity problems.
\section{The algebraic construction}
Let us first mention the following elementary linear algebra fact:
\begin{lemma}\label{lemma:cone}
If $R\in\Sym(\R^N)$ is indefinite, then the cone
\[
 \mathcal C=\{H\in\Sym(\R^N):H[v,v]\ge0\text{ whenever }v^\top Rv=0\}
\]
can equivalently be written as
\[
\mathcal C=\{P+uR:P\ge0,\ u\in\R\}.
\]
\end{lemma}
\begin{proof}
      Let us write $\mathcal K$ to denote the cone $\{P+uR\}$. Clearly $\mathcal K\subseteq \mathcal C$, 
      so it suffices to prove the reverse inclusion. Note that $\mathcal K$ is closed: if $P_j+u_jR$ converges and
$|u_j|\to\infty$, normalisation by $|u_j|$ along a subsequence would
imply that either $R$ or $-R$ is positive semidefinite. 
We claim the dual cones $\mathcal K^*$ and $\mathcal C^*$ satisfy $\mathcal K^*\subseteq \mathcal C^*$, which will prove the lemma.
An elementary calculation proves that $\mathcal K^*=\{B\ge0:\tr(BR)=0\}$ and $\mathcal C^*=\overline{\mathrm{cone}}\{v\otimes v:v^{\top}Rv=0\}$.
For a matrix $B\in\mathcal K^*$, the symmetric matrix $B^{1/2}RB^{1/2}$ is traceless. Therefore it has an orthonormal basis $(e_i)$ with zero diagonal: indeed, 
choosing a unit vector on which the quadratic form vanishes and repeating
on its orthogonal complement provides such a basis. Setting $v_i=B^{1/2}e_i$ then gives
$B=\sum v_i\otimes v_i$ and since $v_i^{\top}Rv_i=0$, it follows that $\mathcal K^*$ is generated by $R$-null rank-one matrices, proving the lemma.
\end{proof}

For $H\in\Sym(\R^4)$, write
$M_{ij}(H)=H_{ii}H_{jj}-H_{ij}^2$ for its principal $2\times2$ minors.
As mentioned in the introduction,
we seek a nonzero, positive linear combination of these minors which
is invariant under addition of $uQ$. There are many such polynomials, but a simple perturbation of $\sigma_2$ suffices:
\begin{lemma}\label{lemma:K}
The quadratic null Lagrangian
\begin{align*}
 K(H)&=\sigma_2(H)+M_{12}(H)+M_{34}(H)\\
 &=2\bigl(M_{12}(H)+M_{34}(H)\bigr)
      +\bigl(M_{13}(H)+M_{14}(H)+M_{23}(H)+M_{24}(H)\bigr),
\end{align*}
defined on $\Sym(\R^4)$,
satisfies $K(H+uQ)=K(H)$. Moreover, if $R\ge aI$ with $a\ge 0$, then $K(R)\ge 8a^2$.
\end{lemma}
\begin{proof}
      Recall that $2\sigma_2(H)=(\tr(H)^2-\tr(H^2))$. Since $\tr Q=0$ and $Q^2=I$, differentiating 
      along $Q$ gives
      \[
      D_Q K(H)=-\tr(HQ)+(H_{11}+H_{22})-(H_{33}+H_{44})=0,
      \]
      proving the first assertion. If $R\ge aI$, each principal
$2\times2$ minor is at least $a^2$ and summing all the minors that make up $K$ we get $8a^2$.
\end{proof}
\begin{remark}
      An alternative construction of $K$, or more precisely a polynomial equivalent for our purposes to $K$, 
      is as follows: take $V$ a two-plane over which $Q$ vanishes. Then $K_V(H):=\det(H_{|V})$, taken 
      with respect to some orthonormal basis of $V$, is a null-lagrangian positive quadratic polynomial 
      satisfying all of our requirements. Our $K$ is simply the average (up to a multiplicative constant)
      of those $K_V$ when varying $V$. As $K$ has the particularly simple form
      $\sigma_2+M_{12}+M_{34}$, this enables us to do the 
      calculations required in the following proofs more explicitly, but our choice is merely one of convenience.\par
      Note that it is the existence of this polynomial $K$ that fails in the $\Sym(2)$ case: no nonzero
      positive quadratic Hessian null lagrangian $K$ that is invariant under addition of multiples of the 
      determinant form $Q$ exists.
\end{remark}
The following
Hessian-minor estimate will play the role of similar $\sigma_k$ estimates from $k$-convexity (see \cite[Section 3]{TW}).

\begin{lemma}\label{lemma:K-bound}
If $f\in C^\infty(B_3)$ is rank-one convex then
\[
 \int_{B_2}K(D^2 f)\lesssim \Lip(f,B_3)^2.
\]
If $f$ is only assumed to be Lipschitz instead of $C^{\infty}$, the same bound holds
on the intersection of $B_2$ with any open set on which $f$ is smooth.
\end{lemma}

\begin{proof}
Assume first $f$ to be smooth and, by scaling, that $\Lip(f,B_3)=1$.
For each minor $M_{ij}$ with $i\ne j$ we have
\[
M_{ij}(f)=f_{ii}f_{jj}-f_{ij}^2=\partial_{ij}(f_i f_j)-\frac12\partial_{ii}(f_j^2)-\frac12\partial_{jj}(f_i^2).
\]
Testing this identity against a nonnegative cutoff $\eta\in C_c^\infty(B_3)$
which is equal to one on $B_2$ gives, by positivity of $K$,
\[
 \int_{B_2}K(D^2f)\le\int\eta K(D^2f)
 \le \sum_{i<j} c_{ij}\int \eta M_{ij}(D^2f),
\]
where $c_{i,j}$ is either $1$ or $2$.
Integrating by parts we get
\[
 \int_{B_2}K(D^2f)\le\int\eta K(D^2f)
 \lesssim \|Df\|_{\infty}^2\|D^2\eta\|_{\LL^1}\le C.
\]
The second assertion follows by mollification and Fatou's lemma.
\end{proof}
In our proof, we will consider the sup-convolution $f_\varepsilon$ of
$f$. We record here a lemma that will prove useful when
calculating $K(D^2 f_\varepsilon)$.
\begin{lemma}\label{lemma:order}
Suppose $H=P+uQ$ with $P$ positive-semidefinite, and write $B=I-\varepsilon H\ge\frac12I$,
with $\varepsilon>0$. Then
\[
 K(HB^{-1})\det B\ge\frac12\varepsilon^2u^4.
\]
In particular, if $|H|>5\tr H$, then
\[
 K(HB^{-1})\det B\ge c\varepsilon^2|H|^4.
\]
\end{lemma}
\begin{proof}
      The hypotheses imply $I/2\le B\le I-\varepsilon uQ$, so both $B$ and $I-\varepsilon u Q$ are positive definite
      and $\varepsilon |u|\le 1/2$. Hence
\[
 HB^{-1}=\frac{B^{-1}-I}{\varepsilon}
 \ge\frac{(I-\varepsilon uQ)^{-1}-I}{\varepsilon}
 =aI+bQ,
\]
where $a=\varepsilon u^2 (1-\varepsilon^2u^2)^{-1}\ge\varepsilon u^2$ and $b=u(1-\varepsilon^2u^2)^{-1}$.
By Lemma \ref{lemma:K} it follows that
\[
 K(HB^{-1})=K(HB^{-1}-bQ)\ge 8a^2.
\]
Since $\det B\ge1/16$, the first claim follows.
To prove the second, we use that $|P|\le\tr P=\tr H$ and $|Q|=2$, so
\[
 |H|\le\tr H+2|u|<\frac15|H|+2|u|,
\]
and $|u|>\frac25|H|$.
\end{proof}
\section{Proof of Theorem \ref{thm:hess-meas} and Corollary \ref{cor:hess-point}}
To deal with the exceptional set where $u$ is not controlled by $\tr(P)$, we will need a standard contact estimate providing a quantitative
bound on the points at which an upper supporting paraboloid of prescribed
curvature is unavailable. This kind of bound is well known (and in particular our estimate follows from \cite[Theorem 1.2]{WZ}): for example, a similar bound is used in Hirsch's proof \cite[Equation 2.1]{Hirsch}, though with a weaker $t^{-\varepsilon}$ estimate\footnote{As we will see later in the proof of
Theorem \ref{thm:hess-meas}, a $t^{-\varepsilon}$ estimate would not suffice to yield integrability, though
the exponent $1$ is not necessary: anything above $2/3$ would do.}
We include a proof for completeness, but we postpone it to after the proof of Theorem \ref{thm:hess-meas}.
\begin{lemma}\label{lemma:contact}
Let $f\in C^\infty(B_8)$, where $B_8\subset\R^4$, with
$\Delta f\ge0$ and $\Lip(f;B_8)\le1$. Define the set $A_t$ by 
\[
A_t=\left\{x\in B_1:
 \sup_{z\in B_4\setminus\{x\}}
 \frac{2[f(z)-f(x)-Df(x)\cdot(z-x)]}{|z-x|^2}\le t\right\}.
\]
Then for $t>1$ we have $|B_1\setminus A_t|\le C/t$.
\end{lemma}
We are now ready to prove the main theorem.
\begin{proof}[Proof of Theorem \ref{thm:hess-meas}]
By a standard slicing argument, we may assume without loss of generality $f$ to be defined on $\R^{2\times2}$; by a smoothing argument we may also assume $f\in C^\infty(B_8)$, $\Lip(f,B_8)=1$.
Lemma \ref{lemma:cone} gives $D^2f=P+uQ$, with $P\ge0$, so $\Delta f\ge0$.
A cutoff gives
\[
 \int_{B_6}\Delta f\le C.
\]
It therefore suffices to estimate $D^2f$ on
\[
 E=\{x\in B_1:|D^2f(x)|>5\Delta f(x)\}.
\]
Take $t\ge 1$, set $\varepsilon=(2t)^{-1}$ and define, for $y\in B_3$, the sup-convolution
\[
 f_\varepsilon(y)=\sup_{|z|\le 1} f(y+z)-\frac{|z|^2}{2\varepsilon}.
\]
As the supremum of a family of rank-one convex $1$-Lipschitz integrands, 
$f_\varepsilon$ is rank-one convex
and $1$-Lipschitz. Let 
$A_t$ be the contact set from Lemma \ref{lemma:contact},
which satisfies
$|B_1\setminus A_t|\le C/t$;
the parabolic upper estimate defining $A_t$ makes the
sup-convolution amenable to a change of variables. More precisely set, for $x\in A_t$,
$T_\varepsilon(x)=x-\varepsilon Df(x)$.
We claim that
$x$ is the unique maximiser in the definition of sup-convolution with $y=T_\varepsilon(x)$.
Since $T_{\varepsilon}(x)$ belongs to $B_2$,
any point $w\in T_\varepsilon(x)+\overline B_1$
lies in $B_4$, and by the definition of $A_t$ we get
\[
 f(w)-\frac{|w-T_\varepsilon(x)|^2}{2\varepsilon}
   -f(x)+\frac{|x-T_\varepsilon(x)|^2}{2\varepsilon}\le-\frac t2|w-x|^2.
\]
Thus $x$ is the unique maximizer, as claimed. It also follows that $T_\varepsilon$ is injective
and indeed locally invertible, since $D^2f\le tI$ implies $DT_{\varepsilon}(x)\ge I/2$.
By construction of $T_\varepsilon$ we have $Df=Df_\varepsilon(T_\varepsilon)$, hence inverting 
in $T_\varepsilon$ and differentiating again yields
\begin{equation}\label{eq:transformed-hessian}
 D^2f_\varepsilon(T_\varepsilon(x))
 =D^2f(x)(I-\varepsilon D^2f(x))^{-1},
\end{equation}
exactly the expression we can estimate via Lemma \ref{lemma:order}.
More precisely, by Lemmas \ref{lemma:order}
and \ref{lemma:K-bound} and the area formula we have
\begin{align}
 c\varepsilon^2\int_{E\cap A_t}|D^2f|^4
 &\le\int_{E\cap A_t}
 K\bigl(D^2f(I-\varepsilon D^2f)^{-1}\bigr)\det(I-\varepsilon D^2f)\notag\\
 &=\int_{T_\varepsilon(E\cap A_t)}K(D^2f_\varepsilon)\le C.
 \label{eq:quartic}
\end{align}
Since $\varepsilon^{-2}=4t^2$, the contact estimate and
\eqref{eq:quartic} imply, for $\lambda\ge1$,
\[
 |E\cap\{|D^2f|>\lambda\}|
 \le |B_1\setminus A_t|
       +|E\cap A_t\cap\{|D^2f|>\lambda\}|\le\frac Ct+\frac{Ct^2}{\lambda^4}.
\]
Choose $t=\lambda^{4/3}$. This yields
$|E\cap\{|D^2f|>\lambda\}|\le C\lambda^{-4/3}$, an integrable tail, and the theorem follows.
\end{proof}
\begin{remark}
      Theorem \ref{thm:hess-meas} can be extended to locally Lipschitz functions $f$ satisfying $\sigma_1(D^2f)\ge 0$ and $K(D^2f)\ge 0$. Indeed, the same proof works, with the exception of Lemma
      \ref{lemma:cone}, though
      one can still obtain Lemma \ref{lemma:order} using the convexity of the cone $\{\tr(H)\ge 0, K(H)\ge 0\}$.
       With a slightly more refined argument, it is also possible to extend it
      to viscosity $(\sigma_1,K)$-convex $\LL^1_{\mathrm{loc}}$ functions.
      Note that $(\sigma_1,K)$-convexity does not imply local Lipschitzness: for example, $\sqrt{|x|}$
      is $(\sigma_1,K)$-convex.
\end{remark}
\begin{proof}[Proof of Lemma \ref{lemma:contact}]
      We borrow $A_t$ and $T_\varepsilon$ from the previous proof and we now take $\varepsilon=t^{-1}$. 
      First observe that on $A_t$ we have, by definition of the set, 
      \[DT_\varepsilon=I-\varepsilon D^2f\ge 0,\]
 and that subharmonicity of $f$ implies
 \[\tr DT_\varepsilon=4-\varepsilon\Delta f\le 4.\]
 It follows that $\det DT_\varepsilon \le 1$.
 We now claim that $B_{1-\varepsilon}\subset T_\varepsilon(A_t)$:
 This, together with the previous estimate, implies the lemma, since
 the area formula gives
\[
 |B_{1-\varepsilon}| \le |T_\varepsilon(A_t)| \le \int_{A_t}
 \det DT_\varepsilon \le|A_t|
\]
and so
\[
 |B_1\setminus A_t| \le |B_1|-|B_{1-\varepsilon}|\le 4|B_1|\varepsilon=\frac{4|B_1|}{t}.
\]
To prove the claim fix $y\in B_{1-\varepsilon}$ and let $z$ maximize $f(w)-|w-y|^2/2\varepsilon$ over $\overline B_4$. Comparing $f(y)$ and $f(z)$ and using Lipschitzness we get
\[
\frac{|z-y|^2}{2\varepsilon}\le f(z)-f(y)\le |z-y|,
\]
so $|z-y|\le 2\varepsilon$; stationarity then gives
$z-y=\varepsilon Df(z)$, hence $T_\varepsilon(z)=y$ and $|z-y|\le\varepsilon$, proving the claim
and thus the lemma. 
\end{proof}
Corollary \ref{cor:hess-point} is a standard consequence of Theorem \ref{thm:hess-meas}:
\begin{proof}[Proof of Corollary \ref{cor:hess-point}]
      If $d=1$ or $m=1$, this is Alexandrov's theorem. Otherwise, Theorem \ref{thm:hess-meas}
gives $Df\in \mathrm{BV}_{\rm loc}$. By the theory of differentiation for $\mathrm{BV}$ functions (see \cite[Theorem 6.1]{EG}) applied to $Df$ we have that, for $R_A$ defined as 
\[R_A(Z)=f(A+Z)-f(A)-Df(A)[Z]-1/2H_f(A)[Z,Z],\]
the estimate
\[
\frac{1}{|B_r|}\int_{B_r}|DR_A|=o(r)
\]
holds
for almost every $A$. Integrating along rays proves
\[
\frac{1}{|B_r|}\int_{B_r}|R_A|=o(r^2).
\]
In other words for almost every $A$ we have
\[
 g_r(Z)\mathrel{:=}\frac{f(A+rZ)-f(A)-rDf(A)(Z)}{r^2}
 \longrightarrow \frac12H_f(A)[Z,Z]
 \quad\text{in }\LL^1_{\rm loc}
\]
as $r$ goes to $0$.
Each $g_r$ is rank-one convex and hence separately convex.
On a compact set $K$, iterating the one-dimensional mean inequality along each separate direction
one obtains a uniform upper bound $M$ for all $g_r$.
On that same set $K$, take $y$ close to $x$ and apply one-dimensional convexity in each direction, bounding 
all vertices other than $x$ with $M$: doing so one obtains
\[
2^N g_r(y)\le g_r(x)+2^N M,
\]
where $N=md$,
and once averaged in a small cube $Q_x$ of admissible $y$ we get, enlarging $M$ so as to bound the local $\LL^1$ norms too,
\[
g_r(x)\ge 2^N \frac{1}{|Q_x|}\int_{Q_x}g_r(y)-2^N M\ge -2^N \frac{M}{|Q_x|}-2^NM.
\]
Choosing $Q_x$ with uniform size yields a uniform lower bound for $g_r$ on $K$.
The upper and lower bounds together imply a uniform bound on the Lipschitz constants, thus $(g_r)$ is locally equicontinuous and therefore locally uniformly converging. Evaluating at
$Z/|Z|$ with $r=|Z|$ proves the claimed expansion.
\end{proof}
\begin{remark}
An application of Alberti's rank-one theorem \cite{Alberti} yields that the singular part of the
Hessian is positive semidefinite. Hence the singular part can only \textit{add} convexity.
\end{remark}
\section{Some more applications}
\subsection{Other cones}
We now turn our attention to two other cones: those generated by quadratic forms, as described in Proposition
\ref{prop:classification} and a cone in $2\times 2$ which is smaller 
than the rank-one cone but nonetheless suffices to ensure $D^2f$ is a measure. We start with the
proof for the quadratic cone:
\begin{proof}[Proof of Proposition \ref{prop:classification}]
      If $\ker q\neq 0$ and $p=0$,this is elementary: choose a linear functional $\ell\neq 0$ which vanishes on 
$\ker q$ and consider $f=g\circ \ell$ where $g:\R\to \R$ is irregular enough. 
The case $\min(p,n)=1$ follows from \cite[Theorem 1.4]{KK} (up to quotienting out the kernel and pulling the counterexample back if necessary), so 
suppose $p,n\ge2$. We may without loss of generality assume
\[
 Q_q= \diag(I_p,-I_n,0_k)
\]
where $Q_q$ is the matrix associated to $q$. Since
the null cone spans $\R^N$, finite $\mathcal D_q$-convex functions
are locally Lipschitz (see \cite[Lemma 2.3]{KK}).
As before, it suffices to prove a uniform local estimate for smooth
functions, so let $f$ be smooth and $\mathcal D_q$-convex on $B_8$, with
$\Lip(f;B_8)\le1$, and write $H=D^2f$.
Lemma \ref{lemma:cone} gives $H=P+uQ_q$, $P\ge0$. Put
\[
 M= \diag(nI_p,pI_n,I_k).
\]
Then $M>0$ and $\tr(MQ_q)=0$, so
\[
 M:H=M:P\ge c_{p,n}|P|\ge0.
\]
Fix the subspace $V$ generated by two positive and two negative
coordinate vectors of $Q_q$ and let $H_{VV}$
denote the restriction of $H$ to $V$. For $u\ne0$, choose a unit vector
$v\in V$ with $Q_qv=\operatorname{sgn}(u)v$. Then
$H[v,v]\ge|u|$, and consequently $|u|\le|H_{VV}|$. Therefore
\begin{equation}\label{eq:q-control}
 |H|\le C_{p,n}\bigl(M:H+|H_{VV}|\bigr).
\end{equation}

The first term has a cutoff bound:
\[
 \int_{B_1}M:H
 \le\int\eta\,M:D^2f
 =-\int M Df\cdot D\eta\le C_{p,n,N},
\]
where $\eta\ge0$ is supported in $B_2$ and equals one on $B_1$.
On each affine $V$-slice, the null cone of $Q_q|_V$ identifies
orthogonally with the rank-one cone in $\R^{2\times2}$.
The four-dimensional estimate and a slicing argument thus give
$\int_{B_1}|H_{VV}|\le C_N$.
Together with \eqref{eq:q-control}, this yields
$\int_{B_1}|D^2f|\le C_{p,n,N}$ and the result follows.
\end{proof}
Something stronger than Theorem \ref{thm:hess-meas} actually holds. In $\R^{2\times 2}$, instead of
requiring full rank-one convexity, it suffices to require convexity on each column and on the "diagonal", where the two columns are the same:
\begin{proposition}
Let $f:\R^{2\times 2}\to \R$ and write it as $f(x,y):\R^2\times \R^2\to \R$.
If $f$ is $\mathcal{D}$-convex where
\[
\mathcal D=\{(u,v):\ u,v\in \R^2\ \text{and either }u=0\text{ or }v=0\text{ or }u=v\}
\] then $f$ has a measure-valued Hessian.
\end{proposition}
The proof is similar to the one of Theorem \ref{thm:hess-meas} so we only sketch it.
\begin{proof}[Sketch]
As in the proof of Theorem \ref{thm:hess-meas}, we only
consider smooth $f$ on $B_8$ with $\Lip(f;B_8)\le1$.
Following the same strategy, we wish to reduce the problem to controlling a scalar component of the Hessian.
This is done by essentially the same linear algebra argument:
write
\[
 H=D^2f=\begin{pmatrix}A&C\\C^\top&B\end{pmatrix}
\]
and $F=A+B+C+C^\top$.
By $\mathcal D$-convexity, $A,B,F\ge0$. Thus
$S:=\tr(A+B+F)\ge0$ has the usual cutoff bound
$\int_{B_1}S\lesssim1$, and controls $A,B$ and
$\Sym C=(F-A-B)/2$.
Writing $C=\Sym C+uJ$, where
$J=(\begin{smallmatrix}0&1\\-1&0\end{smallmatrix})$, gives
$H=uQ_0+R$,
 $Q_0=(\begin{smallmatrix}0&J\\-J&0\end{smallmatrix})$,
and $|R|\lesssim S$, reducing us to having to control $u$ once more.
We hence may restrict ourselves to the set $E:=\{|H|>C_0S\}\cap B_1$, where with $C_0$
sufficiently large we have $H^2\gtrsim |H|^2I$.

The rest of the proof goes along in the same way as that of Theorem \ref{thm:hess-meas}: $\mathcal D$-convexity
implies subharmonicity, so Lemma \ref{lemma:contact} can be applied directly. The sup-convolution argument 
via $T_\varepsilon$ gives
\[
 H_\varepsilon-H
 =\varepsilon H^2(I-\varepsilon H)^{-1}
 \gtrsim\varepsilon|H|^2I,
\]
where $H_\varepsilon=D^2 f_\varepsilon$.
Since $A$ is positive semidefinite, this yields
$M_{12}(H_\varepsilon)\gtrsim\varepsilon^2|H|^4$. The same tail estimates then finish the proof.
\end{proof}
\begin{remark}
    It is interesting to point out that in general one
    cannot hope to have a criterion as simple
    as the one considered in \cite{KK}: indeed, lineality cannot suffice. To see this consider,
    with the same identification of $\R^{2\times 2}$ and $\R^4$, the rank-one cone $\mathcal D$
    and $\mathcal D'=\bigcup_{v\in S}\R v$, where $S=\{e_i\pm e_j:i\in\{1,2\},\ j\in\{3,4\}\}
 \cup\{(1,1,1,1),(1,1,1,-1)\}$.
    Then the two cones of $\mathcal D$ and $\mathcal D'$-convex quadratic forms have the same lineality (and hence the same $\mathcal D$ and $\mathcal D'$-affine functions) and yet by \cite[Theorem 1.4]{KK} there exist $\mathcal D'$-convex functions with unbounded Hessian.
\end{remark}
\subsection{An integral estimate}
Surprisingly, the negative part of the Hessian of a rank-one convex integrand ends up having
better integrability than its positive part. More precisely, following the proof of Theorem
\ref{thm:hess-meas} carefully one gets
\begin{proposition}\label{prop:negative-curvature}
Let $f:\R^{2\times2}\to\R$ be rank-one convex, and write
$D^2f=H_f\mathcal L^4+(D^2f)^s$. Then
\[
 \|H_f^-\|_{\LL^{4/3,\infty}(B_r)}
 \lesssim r^2\Lip(f;B_{8r}),
\]
where $H_f^-=(-H_f)_+$ is the positive semidefinite negative spectral
part.
In particular, $H_f^-\in \LL^q_{\mathrm{loc}}$ for every $1\le q<4/3$.
\end{proposition}
\begin{proof}
Assume without loss of generality that $f$ is smooth and
$\Lip(f;B_8)\le1$. Using the notation of Theorem \ref{thm:hess-meas}, we have $D^2f=P+uQ\ge uQ$ and so it follows that $|H_f^-|\le2|u|$. Thus the first estimate in Lemma \ref{lemma:order} and the
same change of variables give, choosing as before $\varepsilon=(2t)^{-1}$,
\[
 \varepsilon^2\int_{A_t}|H_f^-|^4
 \lesssim\int_{T_\varepsilon(A_t)}K(D^2f_\varepsilon)
 \lesssim1
\]
by Lemma \ref{lemma:K-bound}. Together with Lemma \ref{lemma:contact},
this yields (for $t,\lambda\ge 1$)
\[
 |\{x\in B_1:|H_f^-|>\lambda\}|
 \lesssim t^{-1}+t^2\lambda^{-4}
\]
and the result follows by choosing $t=\lambda^{4/3}$.
\end{proof}
\begin{remark}
As mentioned above, this behaviour is better than that of the positive part of $H_f$. Indeed pointwise
$|H_f(x)|\le |H_f^+(x)|+|H_f^-(x)|$, so considering for example for $p>1$ the integrand $f(A)=|A_{11}|^{2-1/p}$
shows $H_f^+$ cannot in general have better integrability than $\LL^1_{\mathrm{loc}}$. 
\end{remark}
\section{Some more separations between rank-one convexity and convexity}
Now that we know that rank-one convex integrands
have measure-valued Hessians, this opens several natural
follow-up questions related to potential strengthening
of the result. 
We investigate two here: whether every suitably well-behaved rank-one convex integrand can be written locally as the 
difference of convex functions and whether a (suitably well-behaved) rank-one convex integrand $f:L\to \R$
defined on some subspace $L\subset \R^{m\times d}$ can be extended to a rank-one convex integrand on the whole ambient space.
\subsection{Local difference of convex functions}\label{sec:dc}
We start by focusing on the first question, which we now state more precisely:
\begin{center}
Is any rank-one convex integrand locally equal
to the difference of two convex integrands? What if we assume
$f$ to be $\mathrm{C}^{1,\alpha}_{\mathrm{loc}}$?
\end{center}
This turns out not to be the case for any $\alpha<1$:
\begin{theorem}\label{thm:noncvx}
There exists a rank-one convex integrand $f:\R^{2\times2}\to\R$ such that
\[
 f\in\bigcap_{0<\alpha<1}\mathrm{C}^{1,\alpha}_{\mathrm{loc}}(\R^{2\times2})
 \]
and $f$ is not
locally a difference of convex functions.
\end{theorem}
\begin{remark}
In contrast, a $\mathrm C^{1,1}_{\mathrm{loc}}$ rank-one convex integrand is
locally a difference of convex functions. For polyconvex
 integrands, more is true:
 polyconvex integrands are locally semiconvex, and hence can be written locally as the difference of convex functions.
\end{remark}

The construction uses the following extension lemma, extending rank-one convex $C^{1,\omega}$
integrands from $\mathrm{CO}^+(2)$ to $\R^{2\times 2}$,
which will also be useful for the second question. We write
$A=(X,Y)\in\mathrm{CO}^+(2)\oplus \mathrm{CO}^-(2)$ in conformal coordinates.
\begin{lemma}\label{lemma:rough-extension}
Let $\omega$ be a concave Dini modulus and let $g:\mathrm{CO}^+(2)\to \R$ be a $\mathrm{C}^1$ function.
If $Dg$ belongs to  $\mathrm C^{0,\omega}(\mathrm{CO}^+(2))$
then it has a rank-one convex $\mathrm C^1$ extension $f$ and
\[
 f(X,0)=g(X),\qquad Df(X,0)=(Dg(X),0).
\]
If $\omega(t)\lesssim_\alpha t^{\alpha}$ for all $t\in(0,1)$, for all $\alpha\in(0,1)$, then
$f$ has the same regularity as in
Theorem \ref{thm:noncvx}.
\end{lemma}
\begin{proof}
Choose an even smooth mollifier $\eta$ and put
$G(X,r)=g*\eta_r(X)$ and $a(r)=\omega(r)/r$.
Note that
\[
 |G_X-Dg|+|G_r|\lesssim\omega(r)\]
 and
\[
 |G_{XX}|+|G_{Xr}|+|G_{rr}|\lesssim a(r).
\]
Choose $h(0)=h'(0)=0$ and $h''(r)=M(1+a(r))$, and define $f$ as
\[
 f(X,Y)=G(X,|Y|)+h(|Y|)
\]
when $Y\neq 0$, $f(X,0)=g(X)$. By construction, $f$ is  $\mathrm C^1$. Since
$a$ is nonincreasing, $h'(r)/r\ge M(1+a(r))$. For a rank-one
direction $(U,V)$, one has $|U|=|V|$, and hence, outside $\{Y=0\}$,
\[
 D^2f(X,Y)[(U,V),(U,V)]
 \ge\bigl(M/2+(M/2-C)a(|Y|)\bigr)|(U,V)|^2\ge0
\]
provided $M$ is large enough. A nonconstant rank-one line meets $Y=0$
at most once, and its derivative is continuous there; thus $f$ is
rank-one convex. Assume now $\omega$ satisfies the additional
decay condition. Near the $\{Y=0\}$ plane we have
\[|D^2f(X,Y)|\lesssim 1+\frac1{|Y|}\int_0^{|Y|}a(s)\mathrm{d}s,\]
and thus $D^2f\in \LL^p_{\mathrm{loc}}(\R^{2\times 2})$ and the regularity follows 
by standard embedding theorems.
\end{proof}
\begin{proof}[Proof of Theorem \ref{thm:noncvx}]
      Let $v(t)=e^{-1/t}\sin(e^{1/t})\mathbf{1}_{t>0}$, $V(t)=\int_0^t v(s)\mathrm{d}s$ and 
      set $g(X)=V(X_1)$. Since the modulus of continuity of $v$ is $\mathcal O(t\log^2(e/t))$, 
      we can apply Lemma \ref{lemma:rough-extension} to obtain a rank-one convex extension 
      $f$ satisfying the hypothesis of the Theorem. If $f$ were a difference of convex functions, 
      its restriction $V(t)=f((t,0),0)$ would also be, implying (\cite{Hartman}) that $V'(t)$
      would have bounded variation on compact intervals, a contradiction.
\end{proof}
\subsection{Some extension results}\label{sec:ext}
It is well known that a strictly rank-one convex integrand $f:\Sym(2)\to \R$ with uniformly bounded Hessian
admits a rank-one convex extension to $\R^{2\times 2}$ (see \cite{Faraco}, \cite{MS}, \cite{Sverak}, \cite{DalMaso}, \cite{Cagnetti}). Nonetheless, 
as the authors of \cite{Faraco} remark, such a procedure is not known for general rank-one convex integrands and this naturally leads to $Q_3$, which we restate here for the 
comfort of the reader:
\begin{center}
$(Q_3)$\textit{Does a (strictly) rank-one convex integrand 
$f:L\to \R$ admit a rank-one convex extension $\hat f:\R^{m\times d}\to \R$?}
\end{center}
\begin{remark}
Note that this is indeed true if we substitute polyconvexity to rank-one convexity, 
provided $L$ is spanned by its rank-one cone. To see this,
write 
$f=G(M(A))$, where $G$ is convex (and possibly extended-valued) and $M(A)$ is the minor map. Then 
every affine functional on $M(L)$ that is bounded below has to be constant by rank-one affinity of minors, so 
$\mathrm{co}(M(L))=\mathrm{aff(M(L))}$ and a projection onto this affine hull gives a finite polyconvex extension.
\end{remark}

A small modification to Lemma \ref{lemma:rough-extension} yields the following slightly strengthened form of 
the aforementioned extension result: it suffices to control $Df$ via a Dini modulus of continuity.
\begin{corollary}
      Let $\omega$ be a Dini modulus and let $f:\Sym(2)\to \R$ be a rank-one convex function
      such that $Df\in \mathrm{C}^{0,\omega}(\Sym(2))$. Then $f$ admits a rank-one convex extension 
      $F:\R^{2\times 2}\to \R$.
\end{corollary}
\begin{proof}[Proof sketch]
Set $r_j=2^{-j}$ and $f_j=f*\eta_{r_j}$. The mollifier estimates in the proof of the $\mathrm{CO}^+(2)$ extension Lemma \ref{lemma:rough-extension} give $|f_j-f|\le Cr_j\omega(r_j)$ and $|D^2f_j|\lesssim a(r_j)$.

Write $X=\Sym(X)+\Skew(X)$. Then for any rank-one $R$ we have $\det(\Sym R)=-1/2 |\Skew R|^2$. 
Expanding this in an eigen-basis of $\Sym(R)$ and using the rank-one convexity of $f_j$ gives 
\[D^2f_j(\Sym X)[\Sym R,\Sym R]\gtrsim -a(r_j)|\Skew R|^2.\]
Thus, for $K$ large enough,
\[
G_j(X)=f_j(\Sym X)-Cr_j\omega(r_j)+Ka(r_j)|\Skew X|^2
       +M|\Skew X|\sum_{k\ge j}\omega(r_k)
\]
is rank-one convex, where the sum is finite by the assumption that $\omega$ is Dini.
On $|\Skew X|=r_j$ we have $G_{j-1}\ge G_j$, provided $M$ is large enough. Thus we can iteratively define 
$F_0=G_0$, $F_j=\max(F_{j-1}, G_j)\mathbf{1}_{|\Skew X|<r_j}+F_{j-1}\mathbf{1}_{|\Skew X|\ge r_j}$ and this is 
rank-one convex by the boundary estimate. The sequence stabilises off $\Sym(2)$, while on $\Sym(2)$ its limit is $\sup_j(f_j-Cr_j\omega(r_j))=f$.
\end{proof}
On the other hand, we now prove the reason for the lack of a general extension procedure: without global control on $Df$, it is possible 
that no extension exists. Let us first focus on the case of $\R^{2\times 2}$, where the proof is particularly simple:
\begin{proposition}\label{prop:sym}
      There exists a uniformly strictly rank-one convex integrand $f:\Sym(2)\to \R$ which does not admit 
      any finite rank-one convex extension to $\R^{2\times 2}$. The same is true for $\mathrm{Tri}(2)$
      and $ \diag(2)$.
\end{proposition}
\begin{proof}
      Consider the integrand $f(A)\mathrel{:=}\tr^{++}(A)\det^{++}(A)$, where $h^{++}$ is defined as $h$ on the set 
      of matrices with positive eigenvalues and $0$ otherwise. As one easily checks $f$ is rank-one convex
      on all three subspaces considered. To prove 
      it cannot be extended, assume for the sake of contradiction that there exists a rank-one convex integrand $F:\R^{2\times 2}\to \R$ such that $F_{|\Sym(2)}=f$. To obtain a contradiction, we construct a laminate uniformly supported on 
      three points $A_t,B_t,C_t$ such that $C_t$ is diagonal and singular, and their barycentre $D_t$ is diagonal, positive and $\lim_{t\to \infty}f(D_t)=\infty$ while the two limits $\lim_{t\to \infty} A_t$ and $\lim_{t\to \infty}B_t$ are both finite. This will imply a contradiction, since rank-one convexity yields
      \[
      f(D_t)\le \frac{F(A_t)+F(B_t)+F(C_t)}3=\frac{F(A_t)+F(B_t)}3
\]
and taking a limit as $t\to \infty$ proves impossibility. One such choice of matrices is 
\[
A_t=\begin{pmatrix}
      0&-1\\1&\frac5t
\end{pmatrix},\qquad
B_t=\begin{pmatrix}
      0&1\\-1&\frac1t
\end{pmatrix}\qquad
C_t=\begin{pmatrix}
      t&0\\0&0
\end{pmatrix}
\]
and these yield
\[
D_t=\begin{pmatrix}
      \frac t3&0\\ 0&\frac2t
\end{pmatrix}.
\]
Indeed, it is easy to see that $C_t-B_t$ is rank one and so is $C_t+B_t-2A_t$. Thus we can construct the laminate 
by first splitting $D_t$ into $A_t$ with weight 1/3 and $(B_t+C_t)/2$ with weight 2/3 and then evenly split the latter into $B_t$, $C_t$. Thus rank-one convexity of $F$ would yield
\[
\frac23\left(\frac t3+\frac2t\right)=f(D_t)\le \frac{F(A_t)+F(B_t)}3
\]
and letting $t\to \infty$ proves the required contradiction. To obtain a uniformly strict example, it suffices to  consider $\mathrm{tr}^{++}(A)^2\det^{++}(A)+\varepsilon|A|^2$ for $\varepsilon>0$ small enough.
\end{proof}
The argument above can easily be adapted to prove a slightly stronger statement: namely, in the locally uniform topology,
the set of non-extendible rank-one convex integrands 
is comeagre in the set of rank-one convex integrands $f:L\to \mathbb R$ with $L=\Sym(2), \mathrm{Tri}(2), \mathrm{diag}(2)$.

      The integrands considered in Proposition \ref{prop:sym} are also quasiconvex\footnote{
            An integrand $f:L\to \R$ is said to be quasiconvex if for all homogeneous gradient Young measures supported in $L$, $\langle f,\mu\rangle \ge f([\mu])$
      }, so the extension property fails for quasiconvex integrands too. Indeed, quasiconvexity 
      on $\operatorname{Tri}(2)$ and $\diag(2)$ is known to be equivalent to rank-one convexity (see \cite{Kirchheim-Harris} and \cite{Muller}), so those cases are covered automatically. For $\Sym(2)$, the proof of Hessian quasiconvexity 
      is slightly more convoluted: $\det^{++}$
      is known to be $\Sym(2)$-quasiconvex by \cite{Sverak-det}, and we can obtain
      our integrand $\tr^{++}(A)\det^{++}(A)$ from it as follows: given a vector $e$ of norm $1$, a
      positive-definite matrix $H$ and $t\ge 0$ small enough
      \[
\det^{++}(H-te\otimes e)=\det^{++}(H)(1-te^\top H^{-1}e)
\le \det^{++}(H)\left(1-\frac{t}{e^\top He}\right),
      \]
      where the last inequality is an equality if $e$
      is an eigenvector of $H$. Summing over the eigenvectors and optimising in $t$
      yields
      
      \[\underset{\substack{(e_1,e_2)\ \mathrm{orth.}\\ t_1,t_2\ge 0}}{\sup}\sum_i t_i\det^{++}(H-t_ie_i\otimes e_i)=\frac14\tr^{++}(H)\det^{++}(H)
      \]
      and since the supremum and the sum operations preserve Hessian quasiconvexity the result follows.
      In particular, not all Hessian-quasiconvex integrands on $\Sym(2)$ are restrictions of quasiconvex integrands
      on $\R^{2\times 2}$, answering a question posed by the authors of \cite{DalMaso} and proving that 
      some global growth condition on the gradient such as condition $(b)$ in their Theorem 1 is necessary.
      Using the same integrand modifications as in Proposition \ref{prop:sym}, one can also obtain a 
      uniformly strictly hessian-quasiconvex example that is not the restriction of a quasiconvex integrand, and by 
      standard arguments one may even make it smooth.
      The same results apply to higher order $k$-quasiconvexity,
      proving that the hypotheses of \cite{Cagnetti} cannot be removed entirely.
\begin{remark}\label{rem:simplest}
It is in a sense necessary to consider moving laminates: if we do not allow escape at infinity
we cannot get a contradiction, since it is always possible to extend a rank-one convex integrand 
defined on an open neighbourhood of a compact rank-one convex set to the whole ambient space,
see\footnote{There is one technical roadblock to applying the Lemma directly, namely that for a compact $K\subset \Sym(2)$ it is not a priori clear that $K^{rc}$ as a subset of $\R^{2\times 2}$ and as one of $\Sym(2)$ coincide. This is, however, not hard to prove by using the extension result for rank-one convex integrands with uniformly bounded Hessian.} \cite[Lemma 2.3]{MS}. At the same time, a laminate supported on only $2$ points with symmetric 
barycentre also can't work, since it is possible to construct an extension $F$ which satisfies Jensen's inequality on any such laminate, so a moving three-point laminate is the simplest possible construction. To see this, define $F$ as 
\[
F(A)=f(\Sym A)+K(1+|A|^3)|\Skew A|.
\]
Then, given any rank-one split $C=tA+(1-t)B$, let $R$ be the symmetric rank-one matrix 
obtained by removing from $\Sym(A-B)$ the smallest-absolute-value eigenvalue. The intrinsic rank-one split 
$C=t A'+(1-t)B'$, with $A'=C+(1-t)R$ and $B'=C-tR$, has the same barycentre so 
since $|\Sym A-A'|\le |\Skew A|$ and similarly $|B'-\Sym B|\le |\Skew B|$, taking $K$ large enough ensures
$f(A')\le F(A)$ and $f(B')\le F(B)$, implying
\[
f(C)\le t F(A)+(1-t)F(B).
\]
Thus all two-point tests with symmetric barycentre are
compatible with one locally Lipschitz extension.
\end{remark}

We now extend Proposition \ref{prop:sym} to the full Theorem \ref{thm:ext}. The proof will work by reduction to the lower-dimensional cases, which will then be covered by the proposition. For simplicity's sake, we first consider the case of $\mathfrak{sl}(r)$. Here, the case $r=2$ follows directly from Proposition \ref{prop:sym} since $\mathrm{sl}(2)$ is just $\Sym(2)$
up to a rank preserving conjugation;
the general case follows by extending the $\mathfrak{sl}(2)$ example to $\mathfrak{sl}(r)$:
\begin{lemma}\label{lemma:sl}
      Let $r\ge 3$ and let $f:\mathfrak{sl}(2)\to \R$ be a locally semiconvex, rank-one convex integrand whose almost-everywhere
      Hessian has at most quadratic growth. Then $f$ admits a rank-one convex extension $F:\mathfrak{sl}(r)\to \R$.
      If $f$ is uniformly strictly rank-one convex, the extension $F$ can be chosen to have the same property.
\end{lemma}
\begin{proof}
      The hypotheses imply that, for almost every $S$, $D^2f(S)[U,U]\ge -C(1+|S|^2)|\det U|$. Fubini
      and semiconvexity then prove that on every line we have
      \[
      \frac{d^2}{dt^2}f(S+tU)\ge -C(1+|S+tU|^2)|\det U|
      \]
      in distribution. Write now each element of $\mathfrak{sl}(r)$ as 
      \[X=\begin{pmatrix}W&B\\ C&D\end{pmatrix}\] and set
$\tau=\tr D$, $PX=W+\tau I_2/2$ and $Y=(B,C,D,\tau W-BC)$. We claim that, for $K$ large enough, the integrand
$F=f\circ P+K(|Y|^2+|Y|^4)$ is rank-one convex. Differentiating twice along a rank-one line
(note that $Y$ is affine on such a line), it suffices to prove 
\[
(1+|W|^2)|\det(PH)|\lesssim (1+|Y|^2)|Y'|^2.
\]
Since $\det(PH)=-(\tau')^2/4$, this can be rewritten as
\[
(1+|W|^2)|\tau'|^2\lesssim (1+|Y|^2)|Y'|^2,
\]
and to obtain this estimate it suffices to expand $|(\tau W)'|^2$ to get
\[
 (1+|W|^2)|\tau'|^2
 \le |\tau'|^2+|(\tau W)'|^2+|\tau W|\,|(\tau W)''|
 \lesssim (1+|Y|^2)|Y'|^2,
\]
For the strict case, extend $f-\varepsilon|\cdot|^2$ and
add $\varepsilon|X|^2$.
\end{proof}
In particular it follows that $\mathfrak{sl}(r)$ does not have the rank-one extension property and that 
the counterexample can be chosen to be uniformly strictly rank-one convex. $\mathfrak{sl}(r)$ is general enough to describe most possible cases:
\begin{lemma}\label{lemma:sl-reduct}
Let $L\subset \R^{m\times d}$ be spanned by its rank-one matrices and suppose that $L$ is not rectangular.
Then there exists $r\ge 2$ such that there exists a rank-one preserving linear map $T:L\to \mathfrak{sl}(r)$.
Moreover, both $E_{12}$ and $E_{21}$ have a rank-one preimage through $T$.
\end{lemma}
\begin{proof}
    Let $\pi_1(L)=\operatorname{span}\bigcup_{X\in L}\operatorname{im}X$ and
    $\pi_2(L)=\operatorname{span}\bigcup_{X\in L}\operatorname{im}X^\top$.
    We may assume $\R^{m\times d}=\pi_1(L)\otimes \pi_2(L)$.
    Choose a nonzero annihilator $\ell\in L^\perp$ of minimal rank $k$.
    
    If $k=1$, write $\ell=ab^\top$. By construction, any rank-one $X\in L$ satisfies either $a^\top X=0$ or $Xb=0$
    and neither of these maps vanishes identically on $L$. We can thus choose $f$ and $g$ in $L^*\setminus\{0\}$ such that
    $f(X)g(X)=0$ for every rank-one $X\in L$.
    Then $T(X)=f(X)E_{12}+g(X)E_{21}$ satisfies the conclusion.
    
    Suppose now that $\ell$ has rank $k>1$. By minimality, 
    $\ell L^\top\ell \neq \{0\}$, since otherwise 
    every matrix $(\ell a)(b^\top \ell)$ would be a rank-one
    element of $L^\perp$. Since $L$ is spanned by rank-one matrices, there exists one such matrix $X$ such that 
    $\ell X^\top\ell \neq 0$ and by another application of minimality we get a rank-one matrix $Y\in L$ such that 
    $\tr(\ell X^\top\ell Y^\top)\neq 0$.
    Factoring $\ell=QP$ where $Q:\R^k\to \R^m,$
    $P:\R^d\to \R^k$ and setting $T(Z)=PZ^\top Q$
    gives the required map, since $\tr(T(X)T(Y))\neq 0$ and so $T(X)$ and $T(Y)$ can, after rescaling, be simultaneously conjugated to 
    $E_{12}$ and $E_{21}$.
\end{proof}
\begin{proof}[Proof of Theorem \ref{thm:ext}]
      If $L$ is rectangular, the extension follows from the map $X\mapsto P_U XP_V$, so assume 
      $L$ is not rectangular. If $L$ is not spanned by its rank-one cone, choose a linear functional 
      $\ell$ that vanishes on the span. Then $f(X)=\sqrt{|\ell(X)|}+|X|^2$ is uniformly strictly rank-one convex and not Lipschitz, so it cannot
      admit an extension. We may thus assume that $L$ is spanned by its rank-one cone, in which case
      it satisfies the hypotheses of Lemma \ref{lemma:sl-reduct}, so let $T:L\to \mathfrak{sl}(r)$ be the linear 
      map whose existence the lemma implies. Let $h$ be the extension to $\mathfrak{sl}(r)$
      of the uniformly strictly rank-one convex non-extendible integrand on $\mathfrak{sl}(2)$ from Proposition
      \ref{prop:sym} and Lemma \ref{lemma:sl} and set $f(X)=h(TX)$. Let $A=a\otimes b$, $B=c\otimes d$ be rank-one $T$-preimages of $E_{12}$, $E_{21}$. Under the induced identification of the rectangle $\mathrm{span}(a,c)\otimes \mathrm{span}(b,d)$ with $\R^{2\times 2}$, a rank-one convex extension $F$ of $f$ would induce an extension of $h_{|\mathrm{span}(E_{21},E_{12})}$, 
      a contradiction. To obtain a uniformly strictly rank-one convex example, consider 
      $h\circ T+\varepsilon\mathrm{dist}^2(X,\mathrm{span}(A,B))$.
\end{proof}
\section*{Acknowledgements}
The author wishes to thank Lucas Dunckley and Luca Gennaioli for the many insightful conversations on this topic, Daniel Faraco for having shared
some of his notes on the topic, and Jan Kristensen for comments on previous versions
of this work.
\section*{AI disclosure}
No large language models were used in the development of the mathematics in this paper.
The author used LLMs solely for manuscript review and writing suggestions.

\end{document}